\documentclass{amsart}
\usepackage{amssymb}
\usepackage{esint}
\usepackage{filecontents}
\usepackage{mathrsfs} 
\usepackage[all]{xy}
\usepackage[usenames,dvipsnames]{xcolor} 

\usepackage{hyperref} 
\hypersetup{
    colorlinks=true,
    linkcolor=red,
    citecolor=blue,
    }
\usepackage[capitalize, nameinlink]{cleveref}
\numberwithin{equation}{section}
\theoremstyle{plain}
\newtheorem{thm}[equation]{Theorem}
\newtheorem{lem}[equation]{Lemma}

\theoremstyle{definition}

\theoremstyle{remark}

\newtheorem*{rem*}{Remark}
 \DeclareFontFamily{U}{wncy}{}
    \DeclareFontShape{U}{wncy}{m}{n}{<->wncyr10}{}
    \DeclareSymbolFont{mcy}{U}{wncy}{m}{n}
    \DeclareMathSymbol{\Sha}{\mathord}{mcy}{"58}

\DeclareMathOperator{\an}{an} 

\DeclareMathOperator{\Arf}{Arf}
\DeclareMathOperator{\Br}{Br} 
 
\DeclareMathOperator{\Clif}{C} 
\DeclareMathOperator{\ClifInv}{clif} 
 
\DeclareMathOperator{\disc}{disc} 

\DeclareMathOperator{\Gal}{Gal} 

\DeclareMathOperator{\Ind}{ind}

\DeclareMathOperator{\lcm}{lcm}

\DeclareMathOperator{\Mod}{mod}

\DeclareMathOperator{\Per}{per}

\DeclareMathOperator{\Rad}{Rad}
\DeclareMathOperator{\Rank}{rank}
\DeclareMathOperator{\sd}{sd}

\begin{document}
\title[Similarity of quadratic forms]
{Similarity of quadratic forms over global fields in characteristic $ 2 $}
\author{Zhengyao Wu}
\date{\today} 
\address{
Department of Mathematics\\
Shantou University\\
243 Daxue Road\\
Shantou, Guangdong, China 515063}
\email{wuzhengyao@stu.edu.cn}
\subjclass[2010]{Primary: 11E12. Secondary: 11E81, 11E88}
\keywords{similarity, quadratic form, global field, characteristic 2}
\begin{abstract}
Let $ K $ be a global function field of characteristic $ 2 $. 
For each non-trivial place $ v $ of $ K $, let $ K_{v} $ be the completion of $ K $ at $ v $. 
We show that if two non-degenerate quadratic forms are similar over every $ K_{v} $, then they are similar over $ K $. 
This provides an analogue of the version for characteristic not $ 2 $ previously obtained by T.Ono. 
\end{abstract}

\maketitle

\section{Introduction}
Let $ K $ be a global field. 
Let $ \Omega_{K} $  be the set of non-trivial places of $ K $. 
For each $ v\in \Omega_{K} $, let $ K_{v} $ be the completion of $ K $ at $ v $. 
Let $ q $ be a non-degenerate quadratic form over $ K $. 
The Hasse-Minkowski theorem establishes that $ q $ is isotropic over $ K $ if and only if $ q\otimes_{K}K_{v} $ is isotropic over $ K_{v} $ for all $ v\in \Omega_{K} $. For a characteristic not $ 2 $, see, for example, \cite[66:1]{OM}. 
For the characteristic $ 2 $, see \cite[Th.3.2]{Pollak70}. 
As a consequence, two non-degenerate quadratic forms $ f $ and $ g $ over $ K $ are isometric over $ K $ if and only if $ f\otimes_{K}K_{v} $ and $ g\otimes_{K}K_{v} $ are isometric over $ K_{v} $ for all $ v\in \Omega_{K} $. 
We are also interested in the classification of quadratic forms up to similarity. 
Let $ F $ be a field. Two quadratic forms $ f $ and $ g $ over $ F $ are \textit{similar} if there exists $ a\in F^{*} $ such that $ f $ and $ ag $ are isometric. 

Let $ K $ be a global field of characteristic not $ 2 $. 
Two quadratic forms $ f $ and $ g $ are similar over $ K $ if and only if $ f\otimes_{K}K_{v} $ and $ g\otimes_{K}K_{v} $ are similar over $ K_{v} $ for all $ v\in \Omega_{K} $ by  \cite[Th.1]{Ono55}. 
A cohomological version of Ono's theorem is provided in \cite[Th.2.5']{C92}. 
The Hasse principle does not hold for similarity of nonsymmetric bilinear forms by  \cite{C93}. 
The Hasse principle for similarity of hermitian forms fails in general and only holds in special occasions by \cite[Cons.2.11, Th.2.22, Th.3.33]{LUvG}. 
A new proof of Ono's theorem is given by \cite[Prop.8.7]{PR}. 

In this paper, we assume that $ K $ is a global function field of \textit{characteristic $ 2 $} and show that two non-degenerate quadratic forms $ f $ and $ g $ are similar over $ K $ if and only if $ f\otimes_{K}K_{v} $ and $ g\otimes_{K}K_{v} $ are similar over $ K_{v} $ for all $ v\in \Omega_{K} $ (\cref{thm-main}). 
The idea of the proof is the classification of quadratic forms over $ K $ or $ K_{v} $ up to isometry by their rank, discriminant and Clifford invariant, see \cite[Th.2.1]{BFT}. 

Let $ F $ be a field of characteristic $ 2 $ and $ F^{*}=F-\{0\} $ its multiplicative group. 
Let $ H_{2}^{n}(F) $ be the Milne-Kato cohomology as in \cite{Milne76} and \cite{Kato82}. 
We say that $ F $ has \textit{separable dimension} $ \sd_{2}(F)\le 2 $ if $ H_{2}^{3}(L)=0 $ for all finite separable extensions $ L/F $ \cite[p.62, just before Lem.1]{G00}. 
For example, if $ F $ is local or global, then $ \sd_{2}(F)\le2 $.

Let $ V $ be an $ F $-vector space of dimension $ n $. 
Let $ q $ be a non-degenerate quadratic form on $ V $. 
We denote by $ q_{\an} $ the anisotropic part of $ q $, and by $ i_{0}(q) $ the Witt index of $ q $. 
For two non-degenerate quadratic forms $ q $ and $ q' $, we write $ q\simeq q' $ if they are isometric; 
and $ q\sim q' $ if they are similar in $ F $, i.e.~$ q\simeq aq' $ for some $ a\in F^{*} $. 
The \textit{rank} of $ q $ is $ n $. 
Let $ b_{q} $ be the symmetric bilinear form such that $ b_{q}(x,y)=q(x+y)-q(x)-q(y) $ for all $ x,y\in V $. 
Let $ \Rad(b_{q})=\{v\in V~|~b_{q}(v,w)=0 \text{ for all } w\in V\} $ and $  \Rad(q)=\{v\in \Rad(b_{q})~|~q(v)=0\} $.  
We call 
\[
q\text{ \textit{non-degenerate} if }\left \{
\begin{array}{ll}
\dim(\Rad(b_{q}))=1, \Rad(q)=0, & \text{ when }n\text{ is odd}, \\
b_{q}\text{ is non-degenerate}, & \text{ when }n\text{ is even}, \\
\end{array}
\right .
\]
For $ a,b\in F $, let $ [a,b] $ denote the binary quadratic form $ q $ with basis $ e_{1}, e_{2} $ such that $ q(e_{1})=a $, $ q(e_{2})=b $ and $ b_{q}(e_{1}, e_{2})=1 $. 
Let $ \mathbb H=[0,0] $ denote the hyperbolic plane over $ F $. 
By \cite[Cor.7.32]{EKM}, for every non-degenerate quadratic form $ q $ over $ F $,  
\[
q\simeq \left \{
\begin{array}{ll}
\langle d\rangle\perp[a_{1}, b_{1}]\perp\cdots\perp[a_{m}, b_{m}],~d\in F^{*},~ a_{i}, b_{i}\in F,  & \text{ when }n\text{ is odd}, \\
\left [a_{1}, b_{1}\right ]\perp\cdots\perp[a_{m}, b_{m}], ~a_{i}, b_{i}\in F,& \text{ when }n\text{ is even}. \\
\end{array}
\right .
\] 
Let $ \wp\colon F\to F $ be the Artin-Schreier map such that $ \wp(x)=x^{2}+x $ for all $ x\in F $. 
The 
 \textit{discriminant} of $ q $ is defined to be 
\[
\disc(q)=\left \{
\begin{array}{ll}
dF^{* 2}\in F^{*}/F^{* 2},& \text{ when }n\text{ is odd}, \\
\Arf(q)\in F/\wp F,& \text{ when }n\text{ is even}. \\
\end{array}
\right .
\]
To be precise, $ \Arf(\left [a_{1}, b_{1}\right ]\perp\cdots\perp[a_{n}, b_{n}])=a_{1}b_{1}+\cdots+a_{n}b_{n}\Mod\wp F $. 

Let $ \Br(F) $ be Brauer group of $ F $ (see, for example \cite[Sec.2.4]{GS17}) and let $ {_{2}}\Br(F) $ be its $ 2 $-torsion part. 
For a central simple $ F $-algebra $ A $, we denote $ [A] $ its class in $ \Br(F) $;   
its \textit{index} $ \Ind(A)=\Ind([A]) $ is the degree of the unique central division $ F $-algebra $ D $ with $ [D]=[A] $; 
its \textit{period} $ \Per(A)=\Per([A]) $ is the order of $ [A] $ in $\Br(F) $ (see, for example \cite[Sec.2.8]{GS17}). 
The \textit{Clifford invariant} of $ q $ is defined to be
\[
\ClifInv(q)=\left \{
\begin{array}{ll}
\left [ \Clif_{0}(q) \right ]\in {_{2}}\Br(F),& \text{ when }n\text{ is odd}, \\
\left [ \Clif(q) \right ]\in {_{2}}\Br(F),& \text{ when }n\text{ is even}, \\
\end{array}
\right .
\]
where $ \Clif(q) $ is the $ \mathbb Z/2 $-graded Clifford algebra of $ q $ and $ \Clif_{0}(q) $ is the even part of $ \Clif(q) $ (The Clifford invariant $ \ClifInv(q) $ is called the \textit{Witt invariant} of $ q $ in \cite[p.242, para.(-1)]{Knus}). 
The Clifford algebra of the binary form $ [a,b] $ is $ C([a,b])=\displaystyle {a,b\brack F} $, where $ a,b\in F $ and $ \displaystyle{a,b\brack F} $ is the quaternion $ F $-algebra generated by $ i,j $ such that $ i^{2}=a $, $ j^{2}=b $ and $ ij+ji=1 $. 
The \textit{Schur index} of $ q $ is defined to be
\[
\Ind(q)=\left \{
\begin{array}{ll}
\Ind(\Clif_{0}(q)),& \text{ when }n\text{ is odd}, \\
\Ind(\Clif(q)),& \text{ when }n\text{ is even}, \\
\end{array}
\right .
\]

%

\section{The local-global principle}
In this section, we prove the local-global principle for the similarity of non-degenerate quadratic forms over a global function field of characteristic $ 2 $. 

\begin{lem}\label{lem-class}
Let $ F $ be a non-archimedean local field of characteristic $ 2 $. 
Let $ q $ be a non-degenerate quadratic form of rank $ n $ over $ F $. 
Let $ m=\Rank(q|_{\Rad(b_{q})}\perp q_{\an}) $. 

(1) Suppose $ n $ is odd. We have 
\begin{itemize}
\item $ m=1 $ if and only if $ \Ind(q)=1 $. 
\item $ m=3 $ if and only if $ \Ind(q)=2 $.  
\end{itemize}

(2) Suppose $ n $ is even. We have 
\begin{itemize}
\item $ m=0 $ if and only if $ \disc(q)$ is trivial and $ \Ind(q)=1 $. 
\item $ m=2 $ if and only if $ \disc(q)$ is non-trivial and $ \Ind(q)\in \{1,2\} $. 
\item $ m=4 $ if and only if $ \disc(q) $ is trivial and $ \Ind(q)=2 $. 
\end{itemize}
\end{lem}

\begin{proof}
Since $ F $ is a non-archimedean local field of characteristic $ 2 $, every quadratic form of rank $ >4 $ is isotropic \cite[Th.1.1]{Baeza82}. Therefore, $ m\le 4 $. 
Since $ m+2i_{0}(q)=n $, we have $ m\equiv n(\Mod 2) $. 
Since $ F $ is local, we have $ {_{2}}\Br(F)\simeq\mathbb Z/2 $ and $ \Ind(q)=\Ind(\ClifInv(q))=\Per(\ClifInv(q))\in\{1,2\} $. 

(1) When $ n $ is odd, $ \dim(q|_{\Rad(b_{q})})=1$. Suppose $ \disc(q)=dF^{*2} $.

When $ m=1 $, $ q \simeq \langle d\rangle\perp\mathbb H^{(n-1)/2}$, $ d\in F^{*} $. 
From \cite[Lem.2]{MTW}, $ \Clif_{0}(q)\simeq \Clif(d\mathbb H^{(n-1)/2})\simeq \Clif(\mathbb H^{(n-1)/2})$. Thus $ \ClifInv(q)$ is trivial. 

When $ m=3 $, $ q \simeq \langle d\rangle\perp[a,b]\perp\mathbb H^{(n-3)/2}$, $ d\in F^{*} $, $ a,b\in F $ and $ [a,b] $ is anisotropic. 
From \cite[Lem.2]{MTW}, $ \Clif_{0}(q)\simeq \Clif(d[a,b]\perp d\mathbb H^{(n-3)/2})$. We have $d[a,b]\perp d\mathbb H^{(n-1)/2}\simeq [da,d^{-1}b]\perp\mathbb H^{(n-3)/2}$. 
From \cite[just before Prop.5]{MTW}, $ \ClifInv(q)$ is the Brauer class of  $\displaystyle{da,d^{-1}b\brack F}$. The quaternion algebra $\displaystyle{da,d^{-1}b\brack F}$ is not split, otherwise the two quadratic forms $ \langle d\rangle\perp[a,b]\perp\mathbb H^{(n-3)/2}$ and $ \langle d\rangle\perp\mathbb H^{(n-1)/2} $ would have the same rank $ n $, the same discriminant $ dF^{* 2} \in F^{*}/F^{* 2}$ and the same Clifford invariant $ 0 \in {_{2}}\Br(F)$. 
Thus,  by \cite[Th.2.1]{BFT}, $ \langle d\rangle\perp[a,b]\perp\mathbb H^{(n-3)/2}\simeq \langle d\rangle\perp\mathbb H^{(n-1)/2} $, a contradiction to the uniqueness of the Witt decomposition \cite[Th.8.5]{EKM}.

(2) When $ m=0 $, $ q \simeq \mathbb H^{n/2}$ is a hyperbolic space, which has trivial Arf invariant and split Clifford algebra. 

When $ m=2 $, $ q \simeq [a,b]\perp\mathbb H^{(n-2)/2}$, $ a,b\in F $ and $ [a,b] $ is anisotropic. Therefore, $ \Arf(q)=ab $ is nontrival. By \cite[just before Prop.5]{MTW}, $ \ClifInv(q)$ is the Brauer class of $\displaystyle{a,b\brack F}$. Thus $ \Ind(q)\in\{1,2\} $. 

When $ m=4 $, $ q \simeq [a_{1},b_{1}]\perp[a_{2},b_{2}]\perp\mathbb H^{(n-4)/2}$, $ a_{1},b_{1},a_{2},b_{2}\in F $ and $ [a_{1},b_{1}]\perp[a_{2},b_{2}] $ is anisotropic. We have $ \Arf(q)=a_{1}b_{1}+a_{2}b_{2}\Mod \wp F $. 
From \cite[just before Prop.5]{MTW}, $ \ClifInv(q)$ is the Brauer class of $\displaystyle {a_{1},b_{1}\brack F}\otimes {a_{2},b_{2}\brack F} $. 
Similar to the case $ m=3 $, $ \Ind(q)\ne 1 $ and hence $ \Ind(q)=2 $. 
The biquaternion algebra $ \displaystyle{a_{1},b_{1}\brack F}\otimes {a_{2},b_{2}\brack F} $ has Albert form $ \varphi=[1,a_{1}b_{1}+a_{2}b_{2}]\perp[a_{1},b_{1}]\perp[a_{2},b_{2}] $. 
By \cite[(16.5)]{KMRT}, $ \Ind(q)=2 $ if and only if $ i_{0}(\varphi)=1 $. Thus, $ \disc(q)=\Arf(q)=a_{1}b_{1}+a_{2}b_{2} $ is trivial.  
\end{proof}

\begin{lem}\label{lem-similar}
Let $ F $ be a non-archimedean local field of characteristic $ 2 $. 
Let $ V $ be an $ F $-vector space of dimension $ n $. 
Let $ f $ and $ g $ be two non-degenerate quadratic forms on $ V $. 
The two forms $ f\simeq g $ if and only if 
\[\left \{
\begin{array}{ll}
i_{0}(f)=i_{0}(g), &\text{ when }n\text{ is odd}, \\
i_{0}(f)=i_{0}(g) \text{ and }\disc f=\disc g,&\text{ when }n\text{ is even}. \\
\end{array}
\right .\]
\end{lem}

\begin{proof}
If $ f\sim g $, then $ f\simeq ag $ for some $ a\in F^{*} $. Therefore, $ i_{0}(f)=i_{0}(ag)=i_{0}(g) $. When $ n $ is even, $ \Arf(f)=\Arf(ag)=\Arf(g) $ since $ a[b,c]\simeq[ab, a^{-1}c] $ for all $ b,c\in F $. 

Conversely, we suppose $ i_{0}(f)=i_{0}(g) $. Therefore, $n-2i_{0}(f)= n-2i_{0}(g)$ which we denote by $ m $. 

Case 1: $ n $ is odd. Suppose $ \disc f=d_{1}F^{*2} $ and $ \disc g=d_{2}F^{*2} $ for $ d_{1}, d_{2}\in F^{*} $. 

When  $ m=1 $, since $ \langle d_{1}\rangle\sim \langle d_{2}\rangle $, $ f \simeq \langle d_{1}\rangle\perp\mathbb H^{(n-1)/2}\sim\langle d_{2}\rangle\perp\mathbb H^{(n-1)/2}\simeq g$. 

When  $ m=3 $, 
$ \Ind(f)=\Ind(g)=2 $ by \cref{lem-class}(1). 
Since $ F $ is a local field, $ {_{2}}\Br(F) \simeq \mathbb Z/2$. 
Also, $ \Per(\ClifInv(f))=2=\Per(\ClifInv(g)) $. Therefore, $ \ClifInv(f)=\ClifInv(g) $. 
By \cite[Prop.11.4]{EKM}, $ \ClifInv(d_{1}d_{2}^{-1}g)=[\Clif_{0}(d_{1}d_{2}^{-1}g)]=[\Clif_{0}(g)]= \ClifInv(g)$.  The two forms  
$ f $ and $ d_{1}d_{2}^{-1}g $ have the same rank $ n $, the same discriminant $ d_{1}F^{* 2} $ and the same Clifford invariant $ \Clif(g) $.  
Now,  by \cite[Th.2.1]{BFT}, $ f\simeq d_{1}d_{2}^{-1}g $. Thus, $ f\sim g $. 

Case 2: $ n $ is even. 
When $ m=0 $, $ f\simeq \mathbb H^{n/2}\simeq g $ and  hence $ f\sim g $.

When  $ m=2 $, suppose $ f \simeq [a_{1},b_{1}]\perp\mathbb H^{(n-2)/2}$, $ a_{1},b_{1}\in F $ and $ [a_{1},b_{1}] $ is anisotropic; $ g \simeq [a_{2},b_{2}]\perp\mathbb H^{(n-2)/2}$, $ a_{2},b_{2}\in F $ and $ [a_{2},b_{2}] $ is anisotropic. 
We may assume that $ a_{1}\ne 0 $ and $ a_{2}\ne 0 $. 
If $ \Arf(f)=\Arf(g) $, i.e.~$ a_{1}b_{1}=a_{2}b_{2}\Mod \wp F $, then $ [a_{1}, b_{1}]\simeq a_{1}[1, a_{1}b_{1}]\sim a_{2}[1, a_{2}b_{2}]\simeq[a_{2}, b_{2}] $ by \cite[Ex.7.6]{EKM}. 
Thus, $ f \simeq [a_{1},b_{1}]\perp\mathbb H^{(n-1)/2}\sim[a_{2},b_{2}]\perp\mathbb H^{(n-1)/2}\simeq g$ and hence $ f\sim g $. 

When  $ m=4 $, $ \disc(f)=\disc(g) $ is trivial, and we obtain $ \Ind(f)=\Ind(g)=2 $ by \cref{lem-class}(2).  Again, since $ F $ is a local field, $ {_{2}}\Br(F) \simeq \mathbb Z/2$. 
Also, $ \Per(\ClifInv(f))=2=\Per(\ClifInv(g)) $. Therefore, $ \ClifInv(f)=\ClifInv(g) $. 
Now, by \cite[Th.2.1]{BFT}, $ f\simeq g $ and  hence $ f\sim g $. 
\end{proof}


\begin{thm}\label{thm-main}
Let $ K $ be a global function field of characteristic $ 2 $. 
Let $ V $ be a $ K $-vector space of dimension $ n $. 
Let $ f $ and $ g $ be non-degenerate quadratic forms on $ V $. 
Let $ \Omega_{K} $ be the set of all non-trivial places of $ K $. 
For each $ v \in \Omega_{K} $, let $ K_{v} $ be the completion of $ K $ at $ v $ and $ f_{v}, g_{v} $ the scalar extensions of $ f,g $ to $ K_{v} $. 
If $ f_{v}\sim g_{v} $ in $ K_{v} $ for all $ v \in \Omega_{K} $, then $ f\sim g $ in $ K $. 
\end{thm}

\begin{proof}
Since $ K $ is a global function field, every non-trivial place $ v\in \Omega_{K} $ is non-archimedean. 
We are free to use \cref{lem-class} and \cref{lem-similar} for all $ K_{v} $. 

Case 1: $ n $ is odd. 
Suppose $ \disc(f)=d_{1}K^{*2} $ and $ \disc(g)=d_{2}K^{*2}  $ for $ d_{1}, d_{2}\in K^{*} $. 
Let $ a=d_{1}d_{2}^{-1} $. 
Since $ f_{v}\sim g_{v} $ in $ K_{v} $, by \cref{lem-similar}, we have $ i_{0}(f_{v})=i_{0}(g_{v})=i_{0}(ag_{v}) $. 
Also $ \disc(f_{v})=d_{1}K_{v}^{*2}=ad_{2} K_{v}^{*2}=\disc(ag_{v}) $ for all $ v $. 
By \cref{lem-class}(1), $ \Ind(f_{v})=\Ind(ag_{v})\in \{1,2\} $. 
Since $ K_{v} $ is a local field, $ {_{2}}\Br(K_{v}) \simeq \mathbb Z/2$. 
Also, $ \Per(\ClifInv(f_{v}))=\Ind(f_{v})=\Ind(ag_{v})=\Per(\ClifInv(ag_{v})) $. 
Therefore, $ \ClifInv(f_{v})=\ClifInv(ag_{v}) $. 
Now, by \cite[Th.2.1]{BFT}, $ f_{v}\simeq ag_{v} $ in $ K_{v} $ for all $ v $. 
By the Hasse-Minkowski theorem, $ f\simeq ag $. Thus, $ f\sim g $. 

Case 2: $ n $ is even. Since $ f_{v}\sim g_{v} $ in $ K_{v} $ for all $ v $, by \cref{lem-similar}, we have $ i_{0}(f_{v})=i_{0}(g_{v}) $ and $ \disc(f_{v})=\disc(g_{v}) $. 
Let \[ S=\{v \in \Omega_{K}~|~\Rank(f_{v, \an})=\Rank(g_{v,\an})=2\}.\]
If $ v\in \Omega_{K}-S $, i.e.~$ \Rank(f_{v, \an})=\Rank(g_{v,\an})\in \{0,4\} $, then by \cref{lem-class}(2), $ \Ind(f_{v})=\Ind(g_{v})\in \{1,2\} $. 
Since $ K_{v} $ is a local field, $ {_{2}}\Br(K_{v}) \simeq \mathbb Z/2$. 
Also, $ \Per(\ClifInv(f_{v}))=\Ind(f_{v})=\Ind(g_{v})=\Per(\ClifInv(g_{v})) $. 
Therefore, $ \ClifInv(f_{v})=\ClifInv(g_{v}) $. 
By \cite[Th.2.1]{BFT}, $ f_{v}\simeq g_{v} $ for all $ v\in \Omega_{K}-S $. 

Subcase 2a: $ \ClifInv(f_{v})=\ClifInv(g_{v}) $ for all $ v\in S $. 
By  \cite[Th.2.1]{BFT}, $ f_{v}\simeq g_{v} $ for all $ v\in S $. 
We have $ f_{v}\simeq g_{v} $ for all $ v\in \Omega_{K} $. 
By the Hasse-Minkowski theorem, $ f\simeq g $. Hence $ f\sim g $. 

Subcase 2b: $ \Ind(\Clif(f_{s})\otimes \Clif(g_{s}))=2 $ for some $ s\in S $. 
Therefore, $ \Ind(\Clif(f_{s})\otimes \Clif(g_{s}))=2 $. 
Since $ K_{v} $ is local, we have $ {_{2}}\Br(K_{v})\simeq\mathbb Z/2 $ and $ \Ind(\Clif(f_{v})\otimes \Clif(g_{v}))=\Per(\Clif(f_{v})\otimes \Clif(g_{v}))\in \{1,2\} $ for all $ v\in \Omega_{K} $. 
By the Albert-Brauer-Hasse-Noether theorem, $ \Ind(\Clif(f)\otimes \Clif(g))=\lcm\limits_{v\in \Omega_{K}}\Ind(\Clif(f_{v})\otimes \Clif(g_{v}))=2 $. 
Suppose $ Q $ is a quaternion division $ K $-algebra such that $ [Q]=[\Clif(f)\otimes \Clif(g)] $ in $ \Br(K) $. 

Now we take $ d\in \wp^{-1}(\disc(g)) $ and show that $ Q\otimes K(d) $ is split. 
\begin{itemize}
\item For all $ v\in \Omega_{K}-S $, we have $  f_{v} \simeq g_{v} $ from subcase 2a. 
Then $ [Q\otimes K_{v}]=[\Clif(f_{v})\otimes \Clif(g_{v})]=2[\Clif(g_{v})]=0$ in $ \Br(K_{v}) $. 
Therefore, $ [Q\otimes K(d)_{w}] =0$ for all places $ w\in \Omega_{K(d)} $ lying over $ v $. 
\item 
For all $ v\in S $ 
and for all places $ w\in \Omega_{K(d)} $ lying over $ v $, $ K(d)_{w} $ is an extension of $ K_{v} $. We have $ \Rank(f_{K(d)_{w}, \an})\le 2 $ and $ \Rank(g_{K(d)_{w}, \an})\le 2 $. 
Since $  f_{v}\sim g_{v}  $, $ \disc(f_{v})=\disc(g_{v}) $ by \cref{lem-similar}. 
Therefore, $ \disc(f_{K(d)_{w}})=\disc(g_{K(d)_{w}}) $. 
We obtain $ f_{K(d)_{w}}\simeq g_{K(d)_{w}} $ and thus $ [Q\otimes K(d)_{w}]=[\Clif(f_{K(d)_{w}})\otimes \Clif(g_{K(d)_{w}})]=0 $ in $ \Br(K(d)_{w}) $. 

\end{itemize}
We obtain that $ Q\otimes K(d)_{w} $ is split for all $ w\in \Omega_{K(d)} $. 
By the Albert-Brauer-Hasse-Noether theorem, $ Q\otimes K(d) $ is split. 

It follows that $ Q=(K(d), \sigma, a) $ is a cyclic algebra where $ \sigma $ is the generator of $ \Gal(K(d)/K) $ given by $ \sigma(d)=d+1 $ and $ Q= K(d)\oplus cK(d)$ for some $ c\in K^{*} $ with $ c^{2}=a $ and $ dc=cd+c $. Let $ b=a^{-1}\wp(d) $ and consider the binary quadratic form $ [a, b] $, we have $ \disc([a, b])=\disc(g)$ in $ K/\wp K $ and $ \Clif([a, b])=\displaystyle{a, b\brack K}\simeq Q $. 

Next, we show that $ f $ and $ ag $ have the same discriminant. Since $ n $ is even, $ \disc(g_{v})=\disc(ag_{v}) $, by \cref{lem-similar}, $ \disc(f_{v})=\disc(g_{v})=\disc(ag_{v}) $ for all $ v\in \Omega_{K} $. 
Suppose $ q_{1}=[1,\disc(f)] $ and $ q_{2}=[1, \disc(ag)]$. 
We have $ (q_{1})_{v}\simeq (q_{2})_{v} $ for all $ v\in \Omega_{K} $. 
By the Hasse-Minkowski theorem, $q_{1}\simeq q_{2}$ and hence $ \disc(f)=\disc(q_{1})=\disc(q_{2})=\disc(ag) $. 

Now, we show that $ f $ and $ ag $ have the same Clifford invariant. 
Since $ \disc([a, b])=\disc(g)$, $ \disc([a,b]\perp g) $ is trivial. 
By \cite[Prop.5]{MTW}, $\Clif([a,b]\perp g)\simeq M_{2}(\Clif(ag)) $ and hence $ [\Clif([a,b]\perp g)]=[\Clif(ag)] $ in $ \Br(K) $.  
Since the \textit{discriminant module} of $ g $ is trivial (see \cite[Ch.III, (4.2.4) Prop.]{Knus}), we also have $ [\Clif([a,b]\perp g)]=[\Clif([a,b])\otimes \Clif(g)]=[Q\otimes \Clif(g)] $ in $ \Br(K) $ by \cite[Ch.IV, (8.1.1) Prop.~2)]{Knus}. 
We obtain \[ [\Clif(f)]=[\Clif(f)\otimes \Clif(g)\otimes \Clif(g)]=[ Q\otimes \Clif(g)]=[\Clif([a,b]\perp g)]=[\Clif(ag)]. \] 
Hence $\ClifInv(f)=\ClifInv(ag)  $.  

Summarizing subcase 2b, $ f $ and $ ag  $ have the same rank, the same discriminant and the same Clifford invariant; by Theorem 2.1 in \cite{BFT}, $ f\simeq ag $. Hence, $ f\sim g $. 
\end{proof}

%

\section{Acknowledgements}
The author is supported by National Natural Science Foundation of China (No.11701352) and Shantou University Scientific Research Foundation for Talents (No.130-760188). 
The author thanks Yong Hu and Peng Sun for helpful discussions.

\bibliographystyle{plain}
\bibliography{wu}
\end{document}